\documentclass[12pt,a4paper]{article}
\usepackage{amstext}
\usepackage{fancyhdr}
\usepackage{amsmath}

\usepackage{amsfonts,graphicx,bezier, amssymb}
\newenvironment{prf}{\noindent{\bf{Proof.}}~~}{\hfill\rule{1ex}{1ex}\vskip1.5ex}
\newcommand{\Z}{\mathbb Z}

\newcommand{\N}{\mathbb N}

\newcommand{\beqa}{\begin{eqnarray}}
\newcommand{\enqa}{\end{eqnarray}}
\newcommand{\beq}{\begin{eqnarray*}}
\newcommand{\enq}{\end{eqnarray*}}

\newtheorem{qn}{Question}[section]
\newtheorem{rem}{Remark}[section]

\newtheorem{cor}{Corollary}[section]
\newtheorem{propn}{Proposition}[section]
\newtheorem{defn}{Definition}[section]
\newtheorem{exam}{{\bf Example}}[section]
\newtheorem{thm}{Theorem}[section]

\newtheorem{lem}{Lemma}[section]

\newcommand{\noi}{\noindent}

\begin{document}

\begin{center}
 
{\bf A RELATIONSHIP BETWEEN 2-PRIMAL MODULES AND MODULES THAT SATISFY THE RADICAL FORMULA}\\

\vspace*{0.3cm}

David Ssevviiri\\

\vspace*{0.3cm}
Department of Mathematics\\
Makerere University, P.O BOX 7062, Kampala Uganda\\
Email: ssevviiri@cns.mak.ac.ug, ssevviirid@yahoo.com\\

\medskip

\begin{abstract}
 The coincidence of the set of all nilpotent elements of a ring with its prime radical  has a module analogue which
 occurs when  the zero submodule satisfies the radical formula. 
 A ring $R$ is 2-primal if the set of all nilpotent elements of $R$ coincides with its prime radical.
 This fact motivates our study in this paper, namely;   to compare 2-primal submodules and submodules that  satisfy the radical formula.
  A demonstration of the importance of 2-primal   modules in bridging the gap between modules over commutative rings   and modules over
  noncommutative rings is done and new examples of rings and modules that satisfy the radical formula are  also given. 
\end{abstract}
\end{center}
\medskip



\section{Introduction}

\begin{paragraph}\noi
Unless stated otherwise, all rings are unital, associative and not necessarily commutative.
The modules are  left unital.   The set of all positive integers is denoted by $\N$. First, we define
key terms and fix  notation which we later use in the sequel. 
 \end{paragraph}
 
 \begin{paragraph}\noi
A proper ideal $\mathcal{I}$ of a ring $R$ is {\it prime} (resp. {\it completely prime}) if for all ideals $\mathcal{A, B}$ of $R$ 
(resp. $a, b\in R$) $\mathcal{AB}\subseteq \mathcal{I}$ (resp. $ab\in \mathcal{I}$), implies $\mathcal{A}\subseteq \mathcal{I}$ (resp. $a\in \mathcal{I}$)
 or $\mathcal{B}\subseteq \mathcal{I}$ (resp. $b\in \mathcal{I}$). Any completely prime ideal is prime but not conversely; if $R$ is commutative,
  there is no distinction between the two notions. We recall  a generalization of the above two ring theoretic ``primes'' to modules.
  \end{paragraph}

 \begin{defn}\rm A proper submodule $P$ of an $R$-module $M$  for which $RM\not \subseteq P$ is
 \begin{enumerate}
  \item  {\it completely prime} (see \cite{cp}) if  $am\in P$ implies  $m\in P$ or $aM\subseteq P$, for all  $a\in R$ and   $m\in M$;
  \item  {\it prime} (see \cite{Dauns})  if for all  ideals $\mathcal{A}$ of $R$ and  submodules $N$ of $M$,
$\mathcal{A}N\subseteq P$ implies $N\subseteq P$ or $\mathcal{A}M\subseteq P$.
 \end{enumerate}
 \end{defn}
 
\begin{paragraph}\noi  
 Any completely prime submodule is prime but not conversely in general.
If $R$ is commutative, the two notions coincide. A simple module is always prime but it need not be completely
prime. Let $P$ be a submodule of an $R$-module $M$ and $S$ a subset of $M$ such that  $S\not\subseteq P$. By  $(P:S)$   we denote  the set 
$\{r\in R~:~rS\subseteq P\}$. If $P$ is a completely prime submodule of  an $R$-module $M$, then $(P:\{m\})$ is a completely prime ideal
of $R$ such that
 $(P:\{m\})=(P:M)$ for all $m\in M\setminus P$, see \cite[Proposition 2.5]{cp}. On the other hand, if $P$ is a prime submodule of an $R$-module
 $M$, then 
 $(P:\{m\})$ need not be a two sided ideal of $R$ but $(P:N)$ coincides  with $(P:M)$ for all submodules $N$ of $M$ and it is a two sided prime ideal of $R$.
 Evidently, notions of   completely   prime submodules and prime submodules are distinct.
 A module is completely prime (resp. prime)  if its zero submodule is a completely prime (resp. prime) submodule.
 
 \end{paragraph}

 \begin{paragraph}\noi The intersection of all completely
prime (resp. prime) submodules of an $R$-module $M$
 containing the submodule $N$ is called the completely prime (resp. prime) radical of $N$ and is denoted by  $\beta_{co}(N)$ (resp. $\beta(N)$). If $N=0$,
 we  call it the completely prime (resp. prime)  radical of $M$ and write $\beta_{co}(M)$ (resp. $\beta(M)$) instead of $\beta_{co}(0)$ (resp. $\beta(0)$).
 If $M$ has no completely prime (resp. prime) submodules containing a submodule $N$, we write $\beta_{co}(N)=M$ (resp. $\beta(N)=M$).
 \end{paragraph}
 
 \begin{defn}\rm
  A proper submodule $P$ of an   $R$-module $M$ for which $RM\not \subseteq P$ is {\it completely semiprime} (resp. {\it semiprime}) if  $a^2m\in P$
  (resp. $aRam\subseteq P$) implies $am\in P$,   for all $a\in R$ and  $m\in M$. 
 \end{defn}

 \begin{paragraph}\noi 
 A module  is completely semiprime (resp. semiprime)
  if its zero submodule is a completely semiprime (resp. semiprime) submodule. Any completely semiprime submodule is semiprime. The converse does not 
   hold, see \cite[p. 45]{PhDT}.
\end{paragraph}

\subsection{Submodules that satisfy the radical formula}
\begin{paragraph}\noi

 For commutative rings, the set  of all nilpotent elements of a ring $R$ coincides with  the prime radical $\beta(R)$ of $R$ which
  is the intersection of all prime ideals of $R$. In general, if $\mathcal{I}$ is an ideal of a ring $R$ and 
  $$\sqrt{\mathcal{I}}:=\{a\in R~:~a^n\in \mathcal{I}~\text{for some}~n\in \N\},$$ then for any ideal $\mathcal{I}$ of a commutative ring $R$ we have 
  
 \begin{equation}\label{moja}
  \sqrt{\mathcal{I}}=\beta(\mathcal{I}),
 \end{equation}
 where $\beta(\mathcal{I})$ is the intersection of all prime ideals of $R$ containing $\mathcal{I}$.  
  In \cite{MM}, McCasland and Moore  have extended this notion to modules over commutative rings by 
  defining the radical formula of a submodule. The envelope $E_M(N)$ of a submodule $N$ of an $R$-module $M$ is the set
  
  $$E_M(N):=\{rm~:~ r\in R, m\in M~\text{and}~r^km\in N~\text{for some }~k\in \N\}.$$
  
  It is easy to show that if $R$ is a commutative ring and   $M=~_RR$, then $\sqrt{0}=E_M(0)$.
  Since $E_M(N)$ is in general not a submodule of $M$, we consider the submodule $\langle E_M(N)\rangle$  of  $M$ generated by   $E_M(N)$.
  \end{paragraph}
  
  \begin{paragraph}\noi 
  We say that  a submodule  $N$ of an $R$-module   $M$ satisfies the radical formula if 
   
  \begin{equation}\label{mbiri}
   \langle E_M(N)\rangle=\beta(N).
  \end{equation}
     A module satisfies the radical formula if every submodule of $M$ satisfies the radical formula. If every $R$-module satisfies the radical formula, then
  $R$ is also said to satisfy the radical formula. In literature, there has been an intensive study of modules that satisfy the radical formula, see
  \cite{A2008, A2007,  JS, LM, S,   NA, SSM} among others. Unlike commutative rings for which $\sqrt{I}=\beta(I)$
    for any ideal $I$, not all modules over commutative rings satisfy the radical formula.
   \end{paragraph}

   \subsection{2-primal submodules}
\begin{paragraph}\noi 
A  not necessarily commutative  ring $R$ for which  $\sqrt{0}=\beta(R)$ is called a 2-primal ring. This  condition forces $\sqrt{0}$ to be
 an ideal of $R$. It follows from  \cite[Proposition 2.1]{Gary} that a ring  $R$ is 2-primal if and only if 
 $\beta_{co}(R)=\beta(R)$,   where $\beta_{co}(R)$ denotes the completely prime radical of $R$. 
 We remind the reader that $\beta_{co}(R)$   is the  intersection of all completely prime ideals of $R$ and it is called  also the generalized nil
 radical. Similarly, if $\mathcal{I}$ is any ideal of  $R$, then the symbol $\beta_{co}(\mathcal{I})$ stands for 
  the intersection of all completely prime ideals of $R$ containing $\mathcal{I}$. That intersection is called
    the completely prime radical  of   $\mathcal{I}$.  The 2-primal rings were studied by many authors (see, for example, \cite{Gary,  YC, Marks2, Marks}).
 An ideal $\mathcal{I}$ of a ring $R$ is called 2-primal if

 \begin{equation}\label{tatu}
  \beta_{co}(R/\mathcal{I})=\beta(R/\mathcal{I}).
 \end{equation}

In \cite{2p}, a  generalization of 2-primal rings was done   to modules. 
 A submodule $N$ of an $R$-module  $M$ is 2-primal if  
 \begin{equation}
  \beta_{co}(M/N)=\beta(M/N).
 \end{equation} 
  A  module is  2-primal if its zero submodule is 2-primal, i.e., if $\beta_{co}(M)=\beta(M)$. Any module over a commutative ring is 2-primal and 
 a projective  module over a 2-primal ring is 2-primal \cite[Theorem 2.1]{2p}. As 2-primal   rings  bridge the gap between commutative  rings and
 noncommutative rings,
 2-primal modules also bridge the gap between  modules over commutative rings   and  modules over noncommutative rings.
 \end{paragraph}
 
 \subsection{Questions to investigate}
 \begin{paragraph}\noi
 Since a  ring $R$ is 2-primal if and only if $\sqrt{0}=\beta(R)=\beta_{co}(R)$, it is  natural to ask whether a module $M$ is also 2-primal if and only if
 $\langle E_M(0)\rangle =\beta(M)=\beta_{co}(M)$. The answer is  no, all submodules of modules defined over commutative rings are 2-primal but they 
 need not satisfy   the radical formula, i.e., it is possible that $\langle E_M(0)\rangle \not=\beta(M)$ for a 2-primal module $M$. Against this 
  background, we pose the following questions which form the basis of our study in this paper: 
  
  \begin{enumerate}
   \item What is (are) the condition(s) for  a module to be 2-primal if and only if $\langle E_M(0)\rangle =\beta(M)$?
   \item When does a 2-primal submodule satisfy the radical formula?
   \item When does a submodule that satisfies the radical formula become 2-primal?
   \item Whenever  an ideal $I$ of a ring $R$ is 2-primal, the set $\sqrt{I}$ is an ideal of $R$;  when does the set $E_M(N)$ become a
   submodule of $M$ for a given submodule $N$ of $M$?
   \item Can we get modules over noncommutative rings which satisfy the radical formula?
   \item Can we get noncommutative rings which satisfy the radical formula?
  \end{enumerate}
 \end{paragraph}

 \begin{paragraph}\noi 
  Note that, if $N$ is a 2-primal submodule of $M$, $E_M(N)$ is not necessarily a submodule of $M$. Take for instance modules over a commutative ring,
   where each submodule is 2-primal.
 \end{paragraph}

 \begin{paragraph}\noi 
   In Corollary \ref{iff}, we give a necessary and sufficient condition for a module to be 2-primal if and only if $\langle E_M(0)\rangle =\beta(M)$.
  In Propositions \ref{rem} and \ref{pf} which have Lemmas \ref{ddd} and \ref{ccc} respectively as special cases, we give situations for which 
   2-primal submodules satisfy the radical formula. Using these lemmas we are able to obtain modules and rings that satisfy the radical formula 
   (see Theorems \ref{TT2} and \ref{TT1}, respectively). In Corollaries \ref{kkk} and \ref{cd} we give conditions 
  on modules $M$ and their submodules  $N$ for the equality  $E_M(N)=\langle E_M(N)\rangle$.
   
 \end{paragraph}

\section{Main Results}

\begin{lem}\label{l1}
If $N$ is a submodule of an $R$-module $M$, then $$\langle E_M(N)\rangle\subseteq \beta_{co}(N).$$
\end{lem}

\begin{prf}
 Let $m\in E_M(N)$. Then  $m=rn$ for some $r\in R$ and  $n\in M$. Moreover,  there exists  $k\in \N$ such that $r^kn\in N$. So,
 $r^kn\in \beta_{co}(N)$.  Since $\beta_{co}(N)$ is a completely semiprime submodule of $M$, we have $m=rn\in \beta_{co}(N)$. Thus 
 $ E_M(N)\subseteq \beta_{co}(N)$ and finally $\langle E_M(N)\rangle\subseteq \beta_{co}(N).$
\end{prf}

\begin{propn}\label{prop1}
 If $N$ is a completely semiprime submodule of an $R$-module $M$, then $$  E_M(N)=N.$$ 
\end{propn}
\begin{prf}
Obviously, $N\subseteq  E_M(N)$. If  $x\in E_M(N)$, then $x=rm$ and $r^km\in N$ for some $r\in R$, $m\in M$ and $k\in \N$.
As $N$ is completely semiprime we get  $x=rm \in N$.
\end{prf}

\begin{paragraph}\noi 
 In \cite[Proposition 2.1]{AN2012} Azizi and Nikseresht  gave a class of modules $M$ defined over a commutative ring
 for which $E_M(N)$ is always a submodule of $M$. This class consists of all modules $M$
  such that $\beta(N)=N$ for every submodule $N$ of $M$. In Corollary \ref{kkk2} we  give a more general  and bigger class of modules $M$ defined over a not
  necessarily commutative ring   for which  $E_M(N)$ is  a submodule of $M$ for every    submodule $N$ of $M$. The class of modules we provide is that
   of fully completely semiprime modules. It is easy to check that the class of modules $M$  defined over a commutative ring for which 
   $\beta(N)=N$ for each submodule $N$ of $M$ is a class of fully completely semiprime modules since in such a case semiprime is indistinguishable 
   from completely semiprime.      We need Corollary \ref{kkk} first.
\end{paragraph}

\begin{paragraph}\noi 
 Proposition \ref{prop1} implies at once the following:
\end{paragraph}

\begin{cor}\label{kkk}
 For any completely semiprime submodule $N$  of a module $M$, $E_M(N)$ is a submodule of $M$.
\end{cor}

 \begin{cor}\label{kkk2}If all submodules of a module $M$ are completely semiprime, then  
   $E_M(N)$ is  a submodule     of $M$ for any submodule $N$ of $M$.  
 \end{cor}

\begin{cor}\label{cd}
 If $M$ is a  2-primal module, then $E_M(\beta(M))=\beta(M)$. In particular,  $E_M(\beta(M))$ is a submodule of $M$.
\end{cor}
\begin{prf}
As  $M$ is 2-primal, we get  $\beta(M)=\beta_{co}(M)$. Moreover, $\beta_{co}(M)$ is a completely semiprime submodule of $M$
 so the assertion  follows from Proposition \ref{prop1}.
\end{prf}

\begin{cor}\label{c}  If $N$  is a 2-primal submodule of $M$, then
$$\langle E_M(N)\rangle/N\subseteq\langle E_M(\beta(N))\rangle/N =\beta(N)/N.$$  In particular, 
$$\langle E_M(0)\rangle\subseteq\langle E_M(\beta(M))\rangle= \beta(M)$$ for any 2-primal module $M$.
 \end{cor}

\begin{prf}
  Suppose $N$ is a 2-primal submodule of $M$.  Since $\beta_{co}(M/N)$ is a completely semiprime submodule
  of $M/N$ and $\beta(M/N)=\beta_{co}(M/N)$, Proposition \ref{prop1} implies\\ $\langle E_{M/N}(\beta(M/N))\rangle=\beta(M/N)$. But 
  $\langle E_{M/N}(\beta(M/N))\rangle =\langle E_M(\beta(N))\rangle/N$ and 
  $\beta(M/N)=\beta(N)/N$. Hence, $\langle E_M(\beta(N))\rangle/N =\beta(N)/N$.   
  As  $N \subseteq\beta(N)$, we get $\langle E_M(N)\rangle \subseteq\langle E_M(\beta(N))\rangle$ and consequently   
  $\langle E_M(N)\rangle/N\subseteq\langle E_M(\beta(N))\rangle/N$. The second statement follows at once from the 
  first one if we put $N=(0)$. 
 \end{prf}

\begin{paragraph}\noi 
 The next result is a direct consequence of Corollary \ref{c}.
\end{paragraph}

\begin{propn}\label{rem}
 Any 2-primal submodule $N$ of an $R$-module $M$ for which $\beta(N)=N$  satisfies the radical formula. 
 \end{propn}

 \begin{paragraph}\noi 
 Notice that for any 2-primal submodule $N$ of an $R$-module $M$, the conditions:
 $\beta(N)=N$, $\beta_{co}(N)=N$, $\beta_{co}(M/N)=\{\bar{0}\}$ and $\beta(M/N)=\{\bar{0}\}$ are equivalent.
 \end{paragraph}

 \begin{propn} For any $R$-module $M$, the following statements hold:
 \begin{itemize}
  \item[(i)] if $R$ is commutative, then every prime submodule $N$ of $M$ satisfies the radical formula;
  \item[(ii)] a completely prime submodule of $M$ satisfies the radical formula.
 \end{itemize} 
   \end{propn}

   \begin{prf}
    If $R$ is commutative, then prime submodules are completely prime. If a submodule $N$ of $M$ is completely prime, then it is  2-primal and prime. 
   Hence $\beta(N)=N$ and the assertion  follows directly from Proposition \ref{rem}.     
   \end{prf}

\begin{propn}\label{pf}
If $M$ is a 2-primal $R$-module such that $\beta(M)=\beta(R)M$ or $\beta_{co}(M)=\beta_{co}(R)M$, then the zero submodule of $M$ satisfies the
radical formula.
 \end{propn}

\begin{prf}
 Suppose that    $\beta(M)=\beta(R)M$. If $x\in \beta(M)$, then $x=\sum_{i=1}^na_im_i$ with $a_i\in \beta(R)$ and $m_i\in M$.
   Since $\beta(R)$ is nil, each $a_i$ is nilpotent and $a_im_i\in E_M(0)$. Hence, $x\in \langle E_M(0)\rangle$.
   Since $M$ is 2-primal,    Corollary \ref{c} implies $\langle E_M(0)\rangle \subseteq \beta(M)$. A similar 
   proof works if we assume that $\beta_{co}(M)=\beta_{co}(R)M$ .
 \end{prf}

 \begin{exam}\label{hd}\rm
  Projective modules satisfy the equations: $\beta(M)=\beta(R)M$ and  $\beta_{co}(M)=\beta_{co}(R)M$,  see \cite[Proposition 1.1.3]{B}.
 \end{exam}

 \begin{rem}
 If we consider  a module $M$  over a commutative ring, then $\langle E_M(0)\rangle \subseteq \beta(M)$.  We see in Corollary \ref{c} that, this is still
  the case
  when $M$ is 2-primal.    Propositions \ref{rem} and \ref{pf} and Corollary \ref{cd} still hold if 
  we replace  ``$N$ 2-primal'' (resp. ``$M$ 2-primal'')  by 
  ``$R$ is commutative''. This    highlights (together with the results obtained in  \cite{2p}) the importance of 2-primal submodules in bridging
  the gap between modules over commutative rings and modules over noncommutative rings.
 \end{rem}

\begin{paragraph}\noi 
 According to Lee and Zhou in \cite{LZ}, an $R$-module $M$ is reduced if for all $a\in R$ and every $m\in M$, $am=0$ implies $Rm \cap aM=0$. An
 $R$-module is  reduced in  this sense if and only if  for all $a\in R$ and every $m\in M$, $a^2m=0$ implies $aRm=0$ if and only if
 for all $a\in R$ and every $m\in M$, $am=0$
 implies $aRm=0$ and  $a^2m=0$ implies   $am=0$, see  \cite[p.25--26]{PhDT}. This implies that any reduced module in the sense of
 Lee and Zhou is completely semiprime.  
A module $M$  is symmetric if $abm=0$ implies $bam=0$ for $a, b\in R$  and $m\in M$.  An $R$-module $M$ is IFP
(i.e., it has the insertion-of-factor-property)
if whenever $am=0$ for $a\in R$ and 
$m\in M$,  we have $aRm=0$.  An $R$-module $M$ is semi-symmetric if for all $a\in R$ and  every $m\in M$,  $a^2m=0$ implies $(a)^2m=0$ where
 $(a)$ is the ideal of $R$ generated by $a\in R$.
A submodule $N$ of an $R$-module $M$ is Lee-Zhou completely semiprime (resp. symmetric, IFP, semi-symmetric) if in the definition of 
reduced (resp. symmetric, IFP, semi-symmetric) we have $N$ in the place of ``$0$'' and ``$\in$'' or ``$\subseteq$'' (whatever is appropriate)
in the place of ``$=$''. For a detailed account of the origin of symmetric modules, IFP modules and semi-symmetric modules together with their examples,
see \cite{2p}. 
\end{paragraph}

\begin{paragraph}\noi 
The following chart of implications is used in the proof of Lemmas \ref{ccc} and \ref{ddd}; it follows from \cite[Theorems 2.2 and 2.3]{2p}.
For any submodule $P$ of an $R$-module $M$, \\

\begin{table}[h]

\begin{tabular}{ccccccc}
                      &               & $R$ commutative &               &              &            & 2-primal.\cr
                      &               & $\Downarrow$    &               &               &           &  $\Uparrow$\cr

  Lee-Zhou            & $\Rightarrow$ & symmetric       & $\Rightarrow$ & IFP           &$\Rightarrow$   & semi-symmetric \cr
completely semiprime                     &&&&&&\cr
\end{tabular}
\begin{center}
 \begin{tabular}{c}
  Chart 1
 \end{tabular}
\end{center}

 \end{table}

\end{paragraph}

\begin{lem}\label{ccc} For an $R$-module $M$, any one of the following statements implies that the zero submodule of $M$ satisfies the radical formula:
\begin{enumerate}
\item $M$ is 2-primal and free,
\item $M$ is semi-symmetric and free,
 \item $M$ is  semi-symmetric and projective,
  \item $M$ is  IFP and projective,
 \item $M$ is  IFP and free,
 \item $M$ is symmetric and projective,
 \item $M$ is symmetric and free,
 \item $M$  is reduced and projective, 
 \item $M$ is reduced and free,
  \item $R$ is commutative and $M$ is projective,
 \item $R$ is commutative and $M$ is free.
 \end{enumerate}
\end{lem}

\begin{prf}
From the chart of implications above it follows that any of the following implies that $M$ is 2-primal: $R$ is commutative, $M$ is reduced, $M$ is IFP, $M$ is symmetric 
and $M$ is semi-symmetric. Secondly, every free module is projective.  The rest follows from Proposition \ref{pf} and Example \ref{hd}.
\end{prf}

\begin{paragraph}\noi
 Lemma \ref{ccc} recovers \cite[Corollary 8]{JS} which says that a zero submodule of a projective module over a commutative ring satisfies
  the radical formula.
\end{paragraph}

\begin{lem}\label{ddd}
 If  a submodule $N$ of a module $M$ completely semiprime (in the sense of Lee-Zhou),  IFP, symmetric or semi-symmetric such 
 that $\beta(N)=N$, then $N$  satisfies  the radical formula.
\end{lem}

\begin{prf}
This follows from Proposition \ref{rem} and the fact that Lee-Zhou completely semiprime,  IFP, symmetric or semi-symmetric submodules are 2-primal.
\end{prf}

\begin{paragraph}\noi 
 The following lemma  was proved by McCasland and Moore in \cite{MM}. Note that, although they were working with modules over commutative rings, the 
  proof they used still works even when the modules are not defined over a commutative ring.
\end{paragraph}

\begin{lem}{\rm\bf{\cite[Theorem 1.5]{MM}}}\label{th1}
Let $\phi ~:~ M\rightarrow M'$ be an $R$-module epimorphism  and let  $N$ be  a submodule of $M$ such that $N\supseteq \text{Ker}~\phi$.
\begin{enumerate}
 \item[(i)] If $\beta(N)=\langle E_M(N)\rangle$, then $\beta(\phi (N))=\langle E(\phi (N))\rangle$;
 \item[(ii)] If $N'$ is a submodule of $M'$ and $\beta(N')=\langle E(N')\rangle$, then $\beta(\phi ^{-1}(N'))=\langle E(\phi ^{-1}(N'))\rangle$.
\end{enumerate}
\end{lem}

\begin{thm}\label{TT1}
 If the $R$-module $M$ is any one of the modules given in Lemma \ref{ccc} or it is 2-primal and projective, then $M$ satisfies the radical
  formula.
\end{thm}

\begin{prf}
 Let $N$ be a submodule of $M$. For the modules given in Lemma \ref{ccc}, apply Lemma \ref{th1}(ii) and Lemma \ref{ccc} by letting $M'=M/N$ and $N'=N$.
 We know that for a 2-primal and projective module $\beta(M)=\langle E_M(0)\rangle$. When we apply Lemma \ref{th1}(ii)
 by letting $M'=M/N$ and $N'=N$, we get the desired result.
\end{prf}

\begin{paragraph}\noi 
 An alternative proof can be given for the six (6) $R$-modules $M$ in Lemma \ref{ccc} which are free. Recall that every $R$-module $M$ is
 the image of  a free
  $R$-module. This together with Lemma \ref{th1}(i) shows that $M$ satisfies the radical formula.
\end{paragraph}

\begin{cor}
 If $R$ is a semisimple ring such that the $R$-module $M$ is 2-primal, then $M$ satisfies the radical formula.
\end{cor}

\begin{prf}
 If $R$ is semisimple, then the $R$-module $M$ is projective. The rest follows from Theorem \ref{TT1}.
\end{prf}

\begin{cor} 
 If $R$ is a semisimple and commutative ring, then  the $R$-module $M$ satisfies the radical formula. 
\end{cor}

\begin{prf}
 If $R$ is semisimple and commutative, then $M$ is  2-primal and projective and it is sufficient to apply Theorem \ref{TT1}.
\end{prf}

\begin{paragraph}\noi 
 A ring $R$ is absolutely radical if for all $R$-modules $M$, we have $\beta(N)=N$ for each submodule $N$ of $M$.
\end{paragraph}

\begin{thm}\label{TT2}
 If $R$ is an absolutely radical ring such that each submodule  $N$ of the $R$-module $M$ is one of the following:
 Lee-Zhou completely semiprime, IFP, symmetric or semi-symmetric, then $R$ satisfies the radical formula. 
\end{thm}

\begin{prf}
Notice that $R$ is an absolutely radical ring if and only if $\beta(N)=N$ for each submodule $N$ of $M$. The rest follows from Lemma \ref{ddd}.
\end{prf}

\begin{propn}\label{ppp1}
 If a submodule $N$ of an $R$-module $M$ satisfies the radical formula and $\beta_{co}(N)\subseteq \langle E_M(N) \rangle$, then $N$ is 2-primal.
 On the other hand, if a  submodule $N$ of an $R$-module $M$ is 2-primal and $\beta_{co}(N)\subseteq \langle E_M(N) \rangle$, then 
 the zero submodule of the $R$-module  $M/N$ satisfies the radical formula. 
\end{propn}

\begin{prf}
 By hypothesis, $\beta_{co}(N)\subseteq \langle E_M(N) \rangle =\beta(N)$ and in general, $\beta(N)\subseteq \beta_{co}(N)$. 
 It follows that  $\beta_{co}(N) =\beta(N)$ such that $\beta_{co}(M/N)=\beta_{co}(N)/N=\beta(N)/N=\beta(M/N)$.
For the second part, suppose  $\beta(M/N)=\beta_{co}(M/N)$ and  $\beta_{co}(N)\subseteq \langle E_M(N)\rangle$. Then
 $\beta(N)/N=\beta_{co}(N)/N\subseteq \langle E_M(N)\rangle/N$. From Lemma \ref{l1}, $\langle E_M(N)\rangle/N\subseteq \beta_{co}(N)/N.$ This implies 
 $\beta(N)/N=\langle E_M(N)\rangle/N$, i.e., $\beta(M/N)=\langle E_{M/N}(\{\bar{0}\})\rangle$.
\end{prf}

\begin{rem}
 The conditions: (1) $\beta_{co}(N)=N$ (which for example holds
  when $N$ is a completely prime submodule) and  (2) $\langle E_M(N)\rangle=M$ (which for example holds when $M$ is cyclic    and   $R$ is nil or
 $M$ is cyclic and its generator is contained in $N$)
 always guarantee existence of the inclusion $\beta_{co}(N)\subseteq \langle E_M(N) \rangle$.
\end{rem}

\begin{cor}\label{iff}

 The necessary and sufficient condition for the zero submodule  of  an $R$-module $M$ to satisfy the radical formula if and only if $M$ is 2-primal is
 \begin{equation}\label{3}
  \beta_{co}(M)\subseteq \langle E_M(0) \rangle.\footnote{To have $\beta_{co}(M)\subseteq \langle E_M(0) \rangle$ is equivalent to having
  $\beta_{co}(M)= \langle E_M(0) \rangle$ since the reverse inclusion always holds.}
 \end{equation}
\end{cor}

\begin{prf}
 It follows from Proposition \ref{ppp1}.
\end{prf}

   \begin{paragraph}\noi
 The following example  shows that  containment (\ref{3}) in Corollary \ref{iff} does not  hold  in general.
 \end{paragraph}

\begin{exam}\label{E1}\rm 
Define $R=\Z[x]$, $F=R\oplus R$,  $f=(2, x)\in F$ and $\mathcal{P}=2R+Rx$ (which is a maximal ideal of $R$).
If $N=\mathcal{P}f$ and $M=F/N$, then $M$ is completely semiprime    and $\beta(M)=Rf/N\not=0$, see \cite[p. 3600]{JS}.
This shows that $\langle E_M(0)\rangle= 0$
(see Proposition \ref{prop1}) and $\beta_{co}(M)\not=0$ since for modules over a commutative ring,  there is no distinction between completely prime 
(resp. completely semiprime) submodules  and prime (resp. semiprime) submodules.
\end{exam}

\begin{paragraph}\noi 
All submodules of a module defined over a commutative ring are 2-primal but they need not 
  satisfy the radical formula. We do not know of an example of a submodule which
   satisfies  the radical formula  but not 2-primal, although we suspect these examples exist. The motivation of our suspicion is that, 
  for any module $M$, $\beta(M)\subseteq \beta_{co}(M)$ and $\langle E_M(0)\rangle \subseteq \beta_{co}(M)$ and these inclusions are in general
   strict. Hence, it is probably possible that $\beta(M)=\langle E_M(0)\rangle \not\subseteq \beta_{co}(M)$, in which case the zero submodule of 
   $M$ satisfies the radical formula  but not 2-primal.  An affirmative answer to any one of the following questions
 gives us the desired example(s).
\end{paragraph}

\begin{qn}
 Is there a prime module $M$ which is not completely prime  and  $E_M(0)=0$?
\end{qn}

\begin{qn}
 Can we get a completely semiprime module $M$ which is not completely prime and  $\beta(M)=0$?
\end{qn}

 \section*{Acknowledgement}  I thank the referee for  comments that   greatly improved this paper.

\end{document}